\documentclass[preprints,article,accept,moreauthors,pdftex,10pt,a4paper]{Definitions/mdpi} 
\pdfoutput=1

\firstpage{1} 

\history{}

\usepackage{amssymb}
\usepackage{xcolor}

\makeatletter
\def\bstr{b}
\def\bfstr{bf}
\def\cstr{c}
\def\fstr{f}
\def\lst{A,B,C,D,d,E,F,G,H,I,J,K,L,M,N,O,P,Q,R,S,T,U,V,W,X,Y,Z,b}

\newcommand{\MkB}[1]{\expandafter\def\csname\bstr#1\endcsname{\mathbb{#1}}}
  \@for\i:=\lst\do{%
    \expandafter\MkB \i     }
    
\newcommand{\MkBF}[1]{\expandafter\def\csname\bfstr#1\endcsname{\mathbf{#1}}}
  \@for\i:=\lst\do{%
    \expandafter\MkBF \i     }

\newcommand{\MkCal}[1]{\expandafter\def\csname\cstr#1\endcsname{\mathcal{#1}}}
  \@for\i:=\lst\do{%
    \expandafter\MkCal \i     }
    
\newcommand{\MkFrak}[1]{\expandafter\def\csname\fstr#1\endcsname{\mathfrak{#1}}}
  \@for\i:=\lst\do{%
    \expandafter\MkFrak \i     }
    
\makeatother

\usepackage{tensor}

\newcommand{\DNsymb}[4]{\tensor*[^{}_{#1}]{#2}{^{#4}_{#3}}}

\newcommand{\idMat}[1]{\DNsymb{#1}{\hat{1}}{}{}}

\newcommand{\zeroMat}[1]{\DNsymb{#1}{\hat{0}}{}{}}

\newcommand{\DNepsilon}[3]{\DNsymb{#1}{\hat{\epsilon}}{#2}{\,#3}}

\newcommand{\DNeq}[1]{\,\rightsquigarrow_{#1}\,}

\newcommand{\bra}[1]{\left\langle #1\right\vert}
\newcommand{\ket}[1]{\left\vert #1\right\rangle}

\newmuskip\pFqmuskip

\newcommand*\pFq[6][8]{%
  \begingroup 
  \pFqmuskip=#1mu\relax
  \mathchardef\normalcomma=\mathcode`,
  \mathcode`\,=\string"8000
  \begingroup\lccode`\~=`\,
  \lowercase{\endgroup\let~}\pFqcomma
  {}_{#2}F_{#3}{\left[\genfrac..{0pt}{}{#4}{#5};#6\right]}%
  \endgroup
}

\newcommand{\pFqcomma}{{\normalcomma}\mskip\pFqmuskip}

\usepackage{mathtools}

\DeclarePairedDelimiterXPP\seq[2]{}{\big(}{\big)}{_{#2}}{#1}

\Title{Dual Numbers and Operational Umbral Methods}

\Author{Nicolas Behr $^{1,}$*, Giuseppe Dattoli $^{2}$, Ambra Lattanzi $^{2,3}$ and Silvia Licciardi $^{2}$}

\AuthorNames{Nicolas Behr, Giuseppe Dattoli, Ambra Lattanzi, Silvia Licciardi}

\address{%
$^{1}$ \quad Universit\'{e} de Paris, Institut de Recherche en Informatique Fondamentale (IRIF), F-75013 Paris, France; nicolas.behr@irif.fr\\  
$^{2}$ \quad ENEA---Frascati Research Center, Via Enrico Fermi 45, 00044 Rome, Italy; giuseppe.dattoli@enea.it (G.D.), silvia.licciardi@enea.it (S.L.)\\  
$^{3}$ \quad H. Niewodnicza\'{n}ski Institute of Nuclear Physics, Polish Academy of Science, Krak\'{o}w, Poland; ambra.lattanzi@ifj.edu.pl}

\corres{Correspondence: nicolas.behr@irif.fr}

\abstract{Dual numbers and their higher order version are important tools for numerical computations, and in particular for finite difference calculus. Based upon the relevant algebraic rules and matrix realizations of dual numbers, we will present a novel point of view, embedding dual numbers within a formalism reminiscent of operational umbral calculus.}

\keyword{dual numbers, operational methods, umbral image techniques}

\begin{document}

\section{Introduction}

The dual numbers (DNs), introduced during the second half of the XIXth century \cite{clifford, grunwald, segre, yaglom, book_algebra}, can be viewed as abstract entities much like the ordinary complex numbers, and are defined as
\begin{equation}
z=x+\epsilon y\,,\qquad (x,y)\in \mathbb{R}\,
\end{equation}
where the corresponding ``imaginary'' unit or \emph{dual number unit (DNU)} $\epsilon$ is a nilpotent number,
\begin{equation}\label{epsilon_nilpotent}
    \epsilon^2=0\qquad \mbox{and}\qquad\epsilon\neq 0.
\end{equation}
The dual numbers were originally introduced within the context of geometrical studies, and  later exploited to deal with problems in pure and applied mechanics \cite{kotelnikov,study}. %
For instance, it has been demonstrated in~\cite{rooney,hsia,martinez} how to formulate the equations of rigid body motion in terms of just three ``dual'' equations instead of their six ``real'' counterparts (thereby realizing an equivalence between spherical and spatial kinematics). More recently, as further discussed in the present paper, their importance has been recognized in numerical analysis to reduce round-off errors \cite{fike2}. We believe that the use of dual numbers in the applied sciences is not as widespread as it could be, and that many new fields of research would benefit from their relevant introduction. %
An important domain in which they may bring significant novelties is that of the perturbative techniques in classical and quantum mechanics.\\

The main contribution of this paper consists in fixing the underlying algebraic rules of the dual numbers in the wider context of umbral and operational calculus. The paper is organized as follows: %
Section~\ref{sec:two} delivers a basic mathematical introduction to dual numbers. %
Section~\ref{sec:three} is devoted to the description of the computational procedure based upon dual numbers and umbral calculus. %
In Section~\ref{sec:four}, we will provide insight into how this powerful method can be applied to deal with problems arising in different contexts. For illustration, we will consider the Schr\"odinger and the heat equation, cornerstones in the respective fields of Physics. %
Section~\ref{sec:five} provides a conclusion with further considerations for future works.

\section{Higher order dual numbers}\label{sec:two}

The DN algebraic rules~\cite{harkin, ozdemir}, summarized below, are a straightforward consequence of the previous identity~\eqref{epsilon_nilpotent} (with $z=x+\epsilon y$ and $w=u+\epsilon v$):
\begin{equation}\label{operation}
\begin{aligned}
&\mbox{\bf Component-wise algebraic addition}\\
&z+w=x+u+\epsilon(y+\nu)\\
&\mbox{\bf Product}\\
&z\cdot w= xu+\epsilon(x\nu+yu)\\
&\mbox{\bf Inverse}\\
&z^{-1}=\frac{1}{x}\Big(1-\epsilon\,\frac{y}{x}\Big) \quad (x\not=0)\\
&\mbox{\bf Power}\\
&z^{n}=x^n\Big(1+n\epsilon\frac{y}{x}\Big)
\quad (n\in\mathbb{Z}_{\geq0}\,,\; x\neq 0)
\end{aligned}
\end{equation}
While the addition operation is entirely analogous to the component-wise addition operation on two-dimensional vectors, the last three operations (product, inverse and power) characterize the distinguishing special algebraic properties of dual numbers (DNs). The multiplication is commutative, associative and distributive, thus the DNs form a two-dimensional associative and commutative algebra over the real numbers.\\

We will now extend this traditional dual number formalism as motivated by the following type of problem. Consider the Taylor expansion up to some order $k$ (denoted $\approx_k$) of an at least $k$-fold continuously differentiable function $f$ around a point $x$, 
\begin{equation}\label{eq:Ttrunc}
  f(x+y)\approx_k \sum_{m=0}^k \frac{y^m}{m!} f^{(m)}(x)\,.
\end{equation}
Following the \emph{automatic differentiation} paradigm~\cite{rall:1981,rall1996introduction,fike}, since in practice the function $f$ will be implemented in some algorithmic from, it may be advantageous to formulate truncations such as~\eqref{eq:Ttrunc} in terms of \emph{generalized (or higher order) dual numbers}. To this end, let us introduce the families of square matrices $\DNepsilon{k}{\pm}{}$, $\idMat{k}$ and $\zeroMat{k}$ with entries (for $i,j=1,\dotsc,k$)
\begin{equation}\label{eq:DNmatK}
  \left(\DNepsilon{k}{\pm}{}\right)_{i,j}:=\delta_{j,i\pm 1}\,,\qquad
  \left(\idMat{k}\right)_{i,j}:=\delta_{i,j}\,,\qquad
  \left(\zeroMat{k}\right)_{i,j}:=0\,,
\end{equation}
where $\delta_{i,j}$ denotes the Kronecker symbol. It is straightforward to verify that for all $k\geq 2$ and $\ell\geq 0$
\begin{equation}
    \left(\DNepsilon{k}{\pm}{\ell}\right)_{i,j}=\delta_{j,i\pm\ell}\quad \Rightarrow\quad 
           \left(\DNepsilon{k}{\pm}{}\right)^k=0\,.
\end{equation}
Then, under the assumptions~\eqref{eq:Ttrunc}, endowing the function $f(x)$ suitably with a component-wise action on square matrices, we find (for $k\geq 2$)
\begin{equation}\label{eq:TtruncMat}
    f\left(x\, \idMat{k}+y\, \DNepsilon{k}{\pm}{}\right)
    =\sum_{m=0}^{k-1}\frac{1}{m!} y^m f^{(m)}(x)\, \DNepsilon{k}{\pm}{m}\, \,.
\end{equation}
For example, setting $k=2$, reproduces the well-known dual number identity~\cite{fike}
\begin{equation}
    f\left(x\, \idMat{2}+y\, \DNepsilon{2}{+}{}\right)
      =f(x)\,\idMat{2}+y f'(x)\,\DNepsilon{2}{+}{}\\
      =\begin{pmatrix}
        f(x) & y f'(x)\\
        0 & f(x)
        \end{pmatrix}\,.
\end{equation}
It may be verified that e.g. for the choice ``+'' in~\eqref{eq:DNmatK}, the first row of the resulting matrices in~\eqref{eq:TtruncMat} contains the terms of the Taylor expansion up to order $k-1$. More explicitly, introducing the auxiliary notations for the row vector $\bra{{}_k e_1}$ and the column vector $\ket{{}_k\mathbf{1}}$ of length $k\geq 2$,
\begin{equation}
\bra{{}_k e_1}:=(1,0,\dotsc,0)\,,\quad
\ket{{}_k\mathbf{1}}:=
    \begin{pmatrix} 1\\ \vdots\\ 1\end{pmatrix}\,,
\end{equation}
let us define the \emph{order $k$ evaluation operation} acting on some function $F(\DNepsilon{k}{+}{})$ depending on a generalized dual number as
\begin{equation}\label{eq:ev}
\langle F(\DNepsilon{k}{+}{}) \rangle_k := \bra{{}_k e_1} F(\DNepsilon{k}{+}{})\ket{{}_k\mathbf{1}}\,.
\end{equation}
We thus find that
\begin{equation}
 \left\langle f\left(x\, \idMat{k}+y\, \DNepsilon{k}{+}{}\right)\right\rangle_k
  =\sum_{m=0}^{k-1} \frac{y^m}{m!} f^{(m)}(x)\,.
\end{equation}

Recently, expansions such as~\eqref{eq:TtruncMat} have received considerable interest in the field of numerical analysis \cite{berland}. Referring to~\cite{fike} for an overview,  various alternative types of "numbers" have been studied for the purpose of finding optimized numerical schemes for computing $k$-th order derivatives of functions. For example, it has been demonstrated that the use of so-called \emph{hyper-dual numbers} results in first and second derivative calculations that are exact, regardless of the step size~\cite{fike2}.\\

For later convenience, motivated by the identity (for $k\geq 2$)
\begin{equation}
  \exp\left({\DNepsilon{k}{+}{}\, x}\right)
  =\sum_{r=0}^{k-1}\frac{x^r}{r!}\, \DNepsilon{k}{+}{r}\,,
\end{equation}
we may introduce the so-called \emph{truncated exponential polynomials}~\cite{dattoli02} $e_n(x)$ defined through the series
\begin{equation}\label{eq:TEP}
    e_n(x):=\sum_{r=0}^n\frac{x^r}{r!}\,,
\end{equation}
which may be expressed in terms of generalized dual numbers as
\begin{equation}
    \begin{pmatrix}
      e_n(x)\\
      e_{n-1}(x)\\
      \vdots\\
      e_1(x)\\
      1
    \end{pmatrix}
  := \exp\left(\DNepsilon{n+1}{+}{}\, x\right) \ket{{}_{n+1}\mathbf{1}}\,.
\end{equation}
One may thus easily verify the property
\begin{equation}
    e'_n(x)=e_{n-1}(x)\,.
\end{equation}

Having provided a matrix-based extension of ordinary to $k$-th order dual numbers of arbitrary order $k\geq 2$, we will now proceed to develop a computational procedure embedding dual numbers with other techniques inspired by the operational umbral formalism.

\section{Umbral-type methods and Dual Numbers}\label{sec:three}

Starting from this section, we will employ the notational simplification of writing $\epsilon$ for the \emph{dual number unit (DNU)} $\DNepsilon{k}{\pm}{}$ of generalized dual numbers (cf.\ Eq.~\eqref{eq:DNmatK}), making the order $k\geq2$ of the DN explicit only via the analogue of the notation~\eqref{eq:ev}, and masking the matrix nature of $\DNepsilon{k}{\pm}{}$. Thus for some function $F\equiv F(\epsilon)$, we write
\begin{equation}\label{eq:defDNeq}
    F \DNeq{k} G \quad :\Leftrightarrow\quad 
    G=\left(F\big\vert_{\epsilon^{k+1}\to0}\right)\big\vert_{\epsilon\to1}
\end{equation}
for the truncation of $F$ via setting $\epsilon^{k+1}=0$ and afterwards $\epsilon=1$. It is straightforward to verify that this formal definition may be implemented in terms of the matrix representations introduced in Section~\ref{sec:two} via use of~\eqref{eq:ev} as
\begin{equation}
    G=\langle F(\DNepsilon{k+1}{+}{}) \rangle_{k+1}\,.
\end{equation}

Consider then the \emph{dual complex parameter}
\begin{equation}\label{parameter_dual}
    \hat{z}\equiv \hat{z}(a,b):=a+\epsilon b\,.
\end{equation}
Following the principles of umbral calculus\footnote{Albeit the term umbral calculus has been introduced in the seminal papers by Roman and Rota~\cite{roman}, in the following we will make reference to the formalism developed in~\cite{licciardi} which enriches the original formalism with the wealth of techniques derived from the operational calculus.}\cite{roman,licciardi}, we will treat the dual complex parameter $\hat{z}$ as an ordinary algebraic quantity in calculations of integrals, derivatives and other operations, delaying the evaluation of $\hat{z}$ via performing the operation $\DNeq{k}$ to the very end of the computations. We will now illustrate the computational benefits of this approach via a number of examples.

\subsection{Dual shifted Gaussians}
We first consider a Gaussian-type function explicitly containing in its argument the dual complex parameter~\eqref{parameter_dual}, whence the \emph{dual-shifted Gaussian} function
\begin{equation}\label{dual_gaussian}
    f(x)=e^{-\alpha x^2+\hat{z}(a,b)\,x}\,.
\end{equation}
Assuming for instance third order dual numbers (i.e.\ $\epsilon^3=0$), we may write the above function in more conventional terms as
\begin{equation}
     f(x) \DNeq{2} e^{-\alpha x^2+a x}\left[ 1+bx +\tfrac{1}{2} (bx)^2 \right]\,,
\end{equation}
which is easily recognized as the product of a shifted Gaussian with a second degree polynomial.\\

 In full analogy to the umbral operational methods of~ \cite{licciardi}, it is then straightforward to calculate the following integral of the function $f$ of~\eqref{dual_gaussian} via the standard Gaussian integral formula
\begin{equation}\label{integral_dual}
      \int_{-\infty}^{+\infty} f(x)dx=
      \sqrt{\tfrac{\pi}{\alpha}}\,e^{\frac{\hat{z}(a,b)^2}{4\alpha}}=
      \sqrt{\tfrac{\pi}{\alpha}}\,e^{\frac{a^2}{4\alpha}+\frac{ab}{2\alpha}\epsilon+\frac{b^2}{4\alpha}\epsilon^2}\,.
\end{equation}
The term on the right has in fact a definite meaning, since the use of the generating function of the \emph{two variable Hermite polynomials}~\cite{appell}
\begin{subequations}
\begin{align}
    \sum_{n=0}^{\infty}\frac{t^n}{n!}H_n(x,y)&=e^{xt+yt^2}\label{eq:HermiteEGF}\\
    H_n(x,y)&=e^{y\partial_x^2}x^n=n!\sum_{r=0}^{\lfloor\frac{n}{2}\rfloor}\frac{x^{n-2r}y^r}{(n-2r)!r!}\label{eq:defHermite}
    \end{align}
\end{subequations}
permits to cast the r.h.s.\ of~\eqref{integral_dual} into the form
\begin{equation}
    \sqrt{\tfrac{\pi}{\alpha}}\,e^{\frac{a^2}{4\alpha}+\frac{ab}{2\alpha}\epsilon+\frac{b^2}{4\alpha}\epsilon^2}
    =\sqrt{\tfrac{\pi}{\alpha}}\,e^{\frac{a^2}{4\alpha}} 
      \sum_{m\geq 0}\frac{\epsilon^m}{m!}\, H_m\left(\tfrac{ab}{2\alpha},\tfrac{b^2}{4\alpha}\right)
    \DNeq{k} \sqrt{\tfrac{\pi}{\alpha}}\,e^{\frac{a^2}{4\alpha}}  {}_H e_{k}\left(\tfrac{ab}{2\alpha},\tfrac{b^2}{4\alpha}\right)\,.
\end{equation}
Here, ${}_H e_k(x,y)$ denotes the \emph{Hermite-based truncated exponential polynomial}~\cite{dattoli03, dattoli03bis, dattoli04} defined as
\begin{equation}
  {}_H e_k(x,y):=\sum_{r=0}^k\frac{1}{r!}\, H_r\left(x,y\right)\,.
\end{equation}

\subsection{Another form of dual Gaussian}\label{sec:DGv2}

Let us consider as a further example
\begin{equation}\label{eq:dualGaussian}
    g(x):=e^{-\hat{z}(a,b) x^2}
\end{equation}
and the following infinite integral (for $\mbox{Re}(a)>0$)
\begin{equation}\label{infinite_integral}
    \int_{-\infty}^{+\infty}g(x)dx
      =\sqrt{\frac{\pi}{\hat{z}}}=\sqrt{\frac{\pi}{a+\epsilon b}}\\
    \DNeq{k}\sqrt{\frac{\pi}{a}}\sum_{r=0}^k\binom{-\frac{1}{2}}{r}\left(\frac{b}{a}\right)^r\,.
\end{equation}
Here, by invoking the operation $\DNeq{k}$, we obtain a finite series, thus obviating the need to impose any condition on the relevant convergence range. 

\subsection{Examples from symbolic calculus}

The calculus of higher order dual numbers may be further refined via combining it with the wealth of techniques available from the theory of special functions and symbolic calculus as put forward in~\cite{licciardi,nuova1,nuova2, nuova3,nuova4,nuova5}. Consider for illustration the following identity, known from the theory of two-variable Hermite polynomials~\cite{dattoli05},
\begin{equation}\label{eq:HermiteIdA}
    \partial^n_xe^{\alpha x^2}=H_n(2\alpha x, \alpha)e^{\alpha x^2}
\end{equation}
which allows to simplify the task of calculating successive derivatives of the dual Gaussian introduced in~\eqref{eq:dualGaussian}, such as in the computation
\begin{equation}
  \partial^n_xe^{-\hat{z} x^2}
    \overset{\eqref{eq:HermiteIdA}}{=}H_n(-2\hat{z}x,-\hat{z}) e^{-\hat{z}x^2}
    \overset{\eqref{eq:defHermite}}{=}
    n! \sum_{r=0}^{\lfloor \tfrac{n}{2}\rfloor}\sum_{s\geq 0}\frac{(-1)^{n-r+s}2^{n-2r}x^{n-2(r-s)}}{(n-2r)!r!s!} \hat{z}^{n-r+s}.\nonumber 
\end{equation}

Another interesting type of calculus concerns infinite integrals involving rational functions such as
\begin{equation}\label{eq:exC}
 \Phi(x;a,b):=\frac{1}{1+\hat{z}x^2}
    \DNeq{k}
      \frac{1}{1+a x^2}\sum_{r=0}^{k}\Big(-\frac{b x^2}{1+a x^2}\Big)^r\,.
\end{equation}
For example, the infinite integral
\begin{equation}
    \int_{-\infty}^{+\infty}\frac{1}{1+\hat{z}x^2}dx=\frac{\pi}{\sqrt{\hat{z}}}
\end{equation}
may be easily transformed into truncated form in full analogy to the calculation summarized in~\eqref{infinite_integral}.

\subsection{Umbral image type techniques}

Referring to~\cite{behr2018operational} for the precise technical details (compare also~\cite{dattoli05}), suffice it here to provide the following definition for the action of the \emph{formal integration operator} $\hat{\mathbb{I}}$ on the formal variable $v$ (for $\alpha\in \mathbb{C}$):
\begin{equation}
  \hat{\mathbb{I}}(v^{\alpha}):=\frac{1}{\Gamma(\alpha)}\,.
\end{equation}
Then an interesting variant of the example presented in~\eqref{eq:exC} may be obtained as
\begin{equation}
     \hat{\mathbb{I}}\left[\int_{-\infty}^{+\infty}v\Phi(x;a, vb)dx\right]
     =
     \hat{\mathbb{I}}\left[\frac{v\pi}{\sqrt{z(a,v\beta)}}\right]
    \DNeq{k}
      \sqrt{\frac{\pi}{a}}\sum_{r=0}^k \frac{1}{\Gamma(\frac{1}{2}-r)(r!)^2}\left(\frac{b}{a}\right)^r\,.
\end{equation}

In summary, the combination of the concept of higher order dual numbers with techniques from symbolic and umbral-image type calculus appears to offer a large potential in view of novel tools of computation. To corroborate this claim, we will now present some first high-level results in this direction.

\section{Dual numbers and solution of heat- and Schr\"odinger-type equations}\label{sec:four}

Before entering the main topic of this section, let us recall a few useful ``operational rules'',  starting with the \emph{Glaisher identity}~\cite{Dattoli_2008,crofton79}
\begin{equation}\label{eq:Glaisher} 
  e^{\tau \frac{d^2}{dx^2}}e^{-\alpha x^2}=\tfrac{1}{\sqrt{1+4\tau \alpha}}e^{-\frac{\alpha x^2}{1+4\tau \alpha}}\,,
\end{equation}
which can also be understood as the solution of the heat equation with a Gaussian as initial function. It will prove particularly useful in the following to note that according to the definition of the Hermite polynomials $H_n(x,y)$ as given in~\eqref{eq:defHermite}, an alternative interpretation of~\eqref{eq:Glaisher} is provided in terms of the \emph{double lacunary exponential generating function} $\mathcal{H}_{2,0}(\lambda;x,y)$ of the polynomials $H_n(x,y)$, where we employ notations as in~\cite{behr2018explicit}
\begin{equation}\label{eq:HPlacunary}
 e^{\tau \frac{d^2}{dx^2}}e^{-\alpha x^2}=\sum_{n\geq 0}\frac{(-\alpha)^n}{n!}H_{2n}(x,\tau)=\mathcal{H}_{2,0}(-\alpha;x,\tau)\,.
\end{equation}
By specializing eq.~\eqref{eq:Glaisher} to $\alpha=\hat{z}$ (with $\hat{z}=a+\epsilon b$ the dual complex parameter \eqref{parameter_dual}), we obtain the operational identity
\begin{equation}\label{eq:GlaisherDN} 
  e^{\tau \frac{d^2}{dx^2}}e^{-\hat{z} x^2}=\tfrac{1}{\sqrt{1+4\tau \hat{z}}}e^{-\frac{\hat{z} x^2}{1+4\tau \hat{z}}}\,.
\end{equation}
Via the simple factorizations
\begin{equation}
\begin{aligned}
1+4\hat{z}\tau
&=\gamma(a,\tau)\gamma\left(\tfrac{b\epsilon }{\gamma(a,\tau)},\tau\right),\\
\frac{\hat{z}}{1+4\hat{z}\tau}
&=\frac{a}{\gamma(a,\tau)}
+\frac{b\epsilon}{[\gamma(a,\tau)]^2\gamma\left(\tfrac{b\epsilon }{\gamma(a,\tau)},\tau\right)},\\
\gamma(c,\tau)&=1+4c\tau\,,
\end{aligned}
\end{equation}
we may transform the identity~\eqref{eq:GlaisherDN} as 
\begin{equation}
 e^{\tau \frac{d^2}{dx^2}}e^{-\hat{z} x^2}
=\mathcal{H}_{2,0}\left(
-\tfrac{b\epsilon}{[\gamma(a,\tau)]^2};x,\tau\gamma(a,\tau)
 \right)\mathcal{H}_{2,0}(-a;x,\tau)\,.
\end{equation}
By re-inserting the definition of the first double-lacunary EGF, using the Glaisher-identity~\eqref{eq:Glaisher} for the second one and finally truncating to order $k$, we eventually arrive at the compact result
\begin{equation}
e^{\tau \frac{d^2}{dx^2}}e^{-\hat{z} x^2}
\DNeq{k}
\frac{e^{-\frac{\alpha x^2}{\gamma(a,\tau)}}}{\sqrt{\gamma(a,\tau)}}
\sum_{n=0}^k\frac{1}{n!}\left(\tfrac{b}{[\gamma(a,\tau)]^2}\right)^n\,
H_{2n}\left(x,\tau \gamma(a,\tau)\right)\,,\qquad \gamma(a,\tau)=1+4 a \tau\,.
\end{equation}
For example, by evaluating the above expression for second order dual numbers, one finds
\begin{equation}\label{eq:heatSecondOrder}
  e^{\tau \frac{d^2}{dx^2}}e^{-\hat{z} x^2}
  \DNeq{2} \frac{e^{-\frac{a x^2}{\gamma(a,\tau)}}}{\sqrt{\gamma(a,\tau)}}\bigg(1-\frac{b}{\gamma(a,\tau)^2}H_2(x,\tau \gamma(a,\tau))+\frac{b^2}{2\gamma(a,\tau)^4}H_4(x,\tau\gamma(a,\tau))\bigg)\,.
\end{equation}
The above result may be interpreted as the solution of the heat-type equation
\begin{equation}
  \partial_\tau F(x,\tau)=\partial^2_x F(x,\tau)\,,\quad
   F(x,0)=e^{-\hat{z}x^2}\,.
\end{equation}
An analogous problem has been addressed in~\cite{dattoli05} within the framework of a different method. The techniques we have envisaged may be further exploited to treat the \emph{paraxial propagation} of the so-called \emph{flattened distributions}, introduced in~\cite{gori} to study the laser field evolution in optical cavities employing super-Gaussian mirrors~\cite{siegman}. These cavities shape beams whose transverse distribution is not reproduced by a simple Gaussian, but by a function exhibiting a \emph{quasi-constant flat-top}, expressible through a function of the type
\begin{equation}
  \begin{aligned}
    E(x;p):=e^{-|x|^p},\qquad p\in \mathbb{Z}_{>0}\,.
  \end{aligned}
\end{equation}
The paraxial propagation of these beams has less obvious properties than, say, Laguerre or Hermite Gauss modes \cite{siegman}. In order to overcome this drawback, Gori introduced the so-called \emph{flattened beams}~\cite{gori} which permit a fairly natural expansion in terms of Gauss Laguerre/Hermite modes, thus providing a straightforward solution to the corresponding paraxial wave equation.\\

Invoking our formalism as developed so far, we may approximate the aforementioned Gori beams in the form
\begin{equation}\label{eq:Exp}
    E(x;p)\approx Y(x;\alpha|m):=e^{-\alpha x^2}e_m(x^2)\,.
\end{equation}
Here, $e_m(x)$ denotes the truncated exponential polynomials introduced in~\eqref{eq:TEP}, and both parameters $\alpha$ and $m$ depend on $p$ (see~\cite{dattoli05} for further details). Recalling from~\eqref{parameter_dual} the definition $\hat{z}(a,b):=a+b\epsilon$ of the dual complex parameter, the r.h.s.\ of~\eqref{eq:Exp} may be equivalently expressed as
\begin{equation}
    e^{-\hat{z}(\alpha,-1)x^2}\DNeq{m} Y(x;\alpha|m)\,,
\end{equation}
whence as an instance of a \emph{dual Gaussian} as described in Section~\ref{sec:DGv2}. The problem of the relevant propagation can accordingly be reduced to that of an ordinary Gaussian mode, namely to the solution of the Schr\"odinger type equation
\begin{equation}
i\partial_\tau\Psi(x,\tau)=-\partial^2_x\Psi(x,\tau)\,,\quad \Psi(x,0)=Y(x;\alpha|m)\,.
\end{equation}
 Consequently, by invoking the operational identity~\eqref{eq:GlaisherDN}, the paraxial evolution of a flattened beam may be expressed in the form
 \begin{equation}
     \Psi(x,\tau)
      =e^{i\tau\partial^2_x}e^{-\hat{z}(\alpha,-1)x^2}=\frac{1}{\sqrt{1+4i\tau \hat{z}(\alpha,-1)}}
      e^{-\frac{\hat{z}(\alpha,-1)x^2}{1+4i\tau \hat{z}(\alpha,-1)}}\,,
 \end{equation}
 which reproduces indeed the known solution of our problem (compare~\cite{dattoli05}). \\
 
 In a forthcoming paper we will discuss this specific application in further detail by applying the method to the problem of designing super-Gaussian optical systems.

\section{Weyl formula and modified Hermite polynomials}\label{sec:five}

 The wide flexibility of the method we propose is corroborated by the following further example, relevant to the use of operational ordering tools. Let us consider an evolution equation of the form
\begin{equation}\label{eq_evo}
\partial_{\tau}F(x,\tau)=[\gamma\partial_x-\hat{z}x]F(x,\tau)\,,\quad F(x,0)=f(x)\,.
\end{equation}
The relevant procedure for combining  differential calculus with the umbral formalism is described in~~\cite{babusci}. Following this approach, the solution of~\eqref{eq_evo} can be expressed as
\begin{equation}\label{sol_Evo}
    F(x,\tau)=e^{\tau(\gamma\partial_x-\hat{z}x)}f(x)\,.
\end{equation}
In order to evaluate the solution \eqref{sol_Evo} explicitly, we need to suitably ``factorize'' the exponential operator. This so-called \emph{disentanglement operation} may be implemented via the Weyl formula \cite{dattoli06}
\begin{equation}\label{weyl_f}
    e^{\hat{X}+\hat{Y}}=e^{-\frac{1}{2}[\hat{X},\hat{Y}]}e^{\hat{X}}e^{\hat{Y}}\,,
\end{equation}
which is applicable whenever the identities $[\hat{X},[\hat{X},\hat{Y}]]=[[\hat{X},\hat{Y}],\hat{Y}]]=0$ hold.
Applying the Weyl formula~\eqref{weyl_f} to our solution \eqref{sol_Evo}, we obtain
\begin{equation}
     F(x,\tau)=e^{-\frac{1}{2}\tau^2\gamma\hat{z}}e^{-\hat{z}x\tau}f(x+\gamma\tau)\,.  
\end{equation}
Thus the solution at any desired truncation order $k$ may be obtained by invoking the dual number evaluation operation $\DNeq{k}$ of~\eqref{eq:defDNeq}.\\

As already mentioned above, the Weyl formula applies in the example presented because the algebraic structure of the argument of the exponential in~\eqref{sol_Evo} satisfies a special property: the commutators of the associated generators reduce to a constant after the first commutation bracket. A more interesting extension is given by the case in which the generators are embedded into a solvable Lie algebra. In this case, the combined  use of the dual number formalism and of the Wei-Norman ordering method~\cite{Wei} leads to new and interesting results. They deserve a separate treatment that will be reported in a forthcoming paper.\\

As a final example, we define \emph{modified Hermite polynomials} $H_n(x,\hat{z})$, whence ordinary two-variable Hermite polynomials $H_n(x,y)$ as introduced in~\eqref{eq:defHermite} evaluated at $y=\hat{z}$, with $\hat{z}\equiv\hat{z}(a,b)$ the dual complex parameter of~\eqref{parameter_dual}, 
\begin{equation}
   H_n(x,\hat{z})=e^{\hat{z}\partial^2_x} x^n\,.
\end{equation}
It is straightforward to verify that these modified polynomials inherit all the relevant properties from the polynomials $H_n(x,y)$, such as the recurrences
\begin{equation}
  \begin{aligned}
    \partial_x H_n(x,\hat{z})
      &=n H_{n-1}(x,\hat{z})\;,\\
    H_{n+1}(x,\hat{z})
        &=x H_n(x,\hat{z})+2\hat{z}\partial_x H_n(x,\hat{z})\,,
  \end{aligned}
\end{equation}
and we find that they satisfy the second order differential equation
\begin{equation}
    2\hat{z}\partial^2_xH_n(x,\hat{z})+x\partial_xH_n(x,\hat{z})=nH_n(x,\hat{z})\,.
\end{equation}
The explicit form of these truncated polynomials is easily obtained. For example, by using third order dual numbers, which implies
\begin{equation}
e^{\hat{z}\partial^2_x}
  \DNeq{2}
e^{a\partial^2_x}\Big(1+b\partial^2_x+\frac{1}{2}b^2\partial^4_x\Big)\,,
\end{equation}
we find the explicit formula
\begin{equation}
 H_n(x,\hat{z})
  \DNeq{2} H_n(x,a)+b\partial_a H_n(x,a)+\tfrac{1}{2}b^2\partial^2_a H_n(x,a)\,,  
\end{equation} 
where we have invoked the well-known identity
\begin{equation}
  \partial_x^2 H_n(x,y)=\partial_y H_n(x,y)\,.
\end{equation}

\section{Final Comments}

The method we have outlined in this paper offers many computational advantages to treat problems where truncated expansions (not necessarily of Taylor type) of functions are involved.
At its core, the umbral formalism and the notion of higher order dual numbers allow to delay the explicit expansions to later stages in a given calculation, thus opening the possibility to exploit numerous efficient computation strategies from the theory of operational calculus and special functions.\\

The technique we have introduced in this paper is amenable for new applications in various different fields. We have presented herein the solution of parabolic equations in transport problems, and within such a context a fairly important example has been provided by  treating the propagation of flattened beams~\cite{dattoli05,gori} in optics. For brevity, we have just outlined the procedure in terms of a 1-dimensional computation. The relevant extension to the 3-dimensional case does not require any particular conceptual effort, but only a consistent numerical implementation.
In a forthcoming investigation, we will further extend the method and study its potential for treating perturbative problems in classical and quantum mechanics.

\vspace{6pt} 


\authorcontributions{conceptualization, G.D.; methodology, N.B., G.D.; validation, N.B., G.D., S.L.; formal analysis, N.B., G.D., A.L., S.L.; writing -- original draft preparation, N.B., G.D., A.L.; writing -- review and editing, N.B., A.L., S.L.}


\funding{The work of N.B. is supported by funding from the European Union's Horizon 2020 research and innovation programme under the Marie Sk\l{}odowska-Curie grant agreement No~753750. A.L. was supported by the NCN research project OPUS 12 no. UMO-2016/23/B/ST3/01714 and by the NAWA project: Program im. Iwanowskiej PPN/IWA/2018/1/00098. S.L. was supported by a \emph{Enea-Research Center Individual Fellowship}.}

\acknowledgments{N.B. would like to thank the LPTMC (Paris 06) and ENEA Frascati for warm hospitality.}

\conflictsofinterest{The authors declare no conflicts of interest. The founding sponsors had no role in the design of the study; in the collection, analyses, or interpretation of data; in the writing of the manuscript, and in the decision to publish the results.%
}




\reftitle{References}



\end{document}